\documentclass{article}
\usepackage{amssymb}

\newtheorem{theorem}{Theorem}
\newtheorem{acknowledgement}[theorem]{Acknowledgement}

\input{tcilatex}

\begin{document}

\begin{center}
\textbf{A subtly analysis of Wilker inequality}\bigskip

Cristinel Mortici\bigskip

Prof. dr. habil., Valahia University of T\^{a}rgovi\c{s}te, Bd. Unirii 18,
130082

T\^{a}rgovi\c{s}te, Romania, cristinel.mortici@hotmail.com

\[
\]
\end{center}

\textbf{Abstract: }\emph{The aim of this work is to improve Wilker
inequalities near the origin and }$\pi /2.$\bigskip



\section{Introduction and Motivation}

In 1989, J. B. Wilker \cite{w} presented the following inequality for $x\in
\left( 0,\pi /2\right) $%
\[
\left( \frac{\sin x}{x}\right) ^{2}+\frac{\tan x}{x}>2.
\]%
This inequality is of great practical importance and it was extended in
different forms in the recent past. We refer to \cite{chh}-\cite{z3} and all
references therein. As Wilker \cite{w} asked about the largest constant $c$
such that%
\[
\left( \frac{\sin x}{x}\right) ^{2}+\frac{\tan x}{x}>2+cx^{3}\tan x~,\text{
\ \ }x\in \left( 0,\pi /2\right) ,
\]%
Sumner, Jagers, Vowe and Anglesio \cite{su} proved the following sharp
inequality for $x\in \left( 0,\pi /2\right) $%
\begin{equation}
2+\frac{16}{\pi ^{4}}x^{3}\tan x<\left( \frac{\sin x}{x}\right) ^{2}+\frac{%
\tan x}{x}<2+\frac{8}{45}x^{3}\tan x.  \label{a}
\end{equation}%
Constants $8/45$ and $16/\pi ^{4}$ are somehow motivated, since they are the
limits at $0,$ respective $\pi /2$ of the function%
\[
x\mapsto \frac{\left( \frac{\sin x}{x}\right) ^{2}+\frac{\tan x}{x}-2}{%
x^{3}\tan x}.
\]%
Recently, Chen and Cheung \cite{chh} proved that this function decreases
monotonically on $\left( 0,\pi /2\right) $ from $8/45$ to $16/\pi ^{4}.$

It is true that inequalities (\ref{a}) and some of recent improvements are
nice through their symmetric form, but let us remember the practical
importance of an inequality which is to provide some bounds for a given
expression. In case of (\ref{a}) observe that near $\pi /2,$ the right-hand
side inequality becomes weak, since%
\[
\lim_{x\rightarrow \left( \pi /2\right) _{-}}\left( \left( \frac{\sin x}{x}%
\right) ^{2}+\frac{\tan x}{x}-\left( 2+\frac{8}{45}x^{3}\tan x\right)
\right) =-\infty .
\]%
As a consequence, if we are interested in finding good approximations of
expression $\left( \sin ^{2}x\right) /x^{2}+\left( \tan x\right) /x$ in
terms of Wilker inequality, then near zero, the constant $8/45$ should be
used, while near $\pi /2,$ the best choice is $16/\pi ^{4}.$

In other words, accurate approximations near zero are obtained of the form%
\[
\left( \frac{\sin x}{x}\right) ^{2}+\frac{\tan x}{x}\approx 2+\left( \frac{8%
}{45}+\lambda \left( x\right) \right) x^{3}\tan x,
\]%
with $\lambda \left( x\right) \rightarrow 0,$ as $x\rightarrow 0,$ while
accurate approximations near $\pi /2$ are obtained of the form%
\[
\left( \frac{\sin x}{x}\right) ^{2}+\frac{\tan x}{x}\approx 2+\left( \frac{16%
}{\pi ^{4}}+\mu \left( x\right) \right) x^{3}\tan x,
\]%
with $\mu \left( x\right) \rightarrow 0,$ as $x\rightarrow \pi /2,$ $x<\pi
/2.$

\section{The Results}

In the light of the discussion from the previous section, we propose the
following new results.\bigskip

\textbf{Theorem 1. }\emph{For every }$x\in \left( 0,1\right) ,$ \emph{we have%
}%
\begin{equation}
2+\left( \frac{8}{45}-a\left( x\right) \right) x^{3}\tan x<\left( \frac{\sin
x}{x}\right) ^{2}+\frac{\tan x}{x}<2+\left( \frac{8}{45}-b\left( x\right)
\right) x^{3}\tan x,  \label{we}
\end{equation}%
\emph{where}%
\[
a\left( x\right) =\frac{8}{945}x^{2}\ ,\ \ \ b\left( x\right) =\frac{8}{945}%
x^{2}-\frac{16}{14175}x^{4}.
\]

\textbf{Theorem 2. }\emph{For every }$x\in \left( \frac{\pi }{2}-\frac{1}{3},%
\frac{\pi }{2}\right) $ \emph{in the left-hand side}$,$ \emph{and for every }%
$x\in \left( \frac{\pi }{2}-\frac{1}{2},\frac{\pi }{2}\right) $\emph{\ in
the right-hand side, the following inequalities hold true:}%
\begin{equation}
2+\left( \frac{16}{\pi ^{4}}+c\left( x\right) \right) x^{3}\tan x<\left( 
\frac{\sin x}{x}\right) ^{2}+\frac{\tan x}{x}<2+\left( \frac{16}{\pi ^{4}}%
+d\left( x\right) \right) x^{3}\tan x  \label{bb}
\end{equation}%
\emph{where}%
\[
c\left( x\right) =\left( \frac{160}{\pi ^{5}}-\frac{16}{\pi ^{3}}\right)
\left( \frac{\pi }{2}-x\right) ,\ \ \ d\left( x\right) =c\left( x\right)
+\left( \frac{960}{\pi ^{6}}-\frac{96}{\pi ^{4}}\right) \left( \frac{\pi }{2}%
-x\right) ^{2}.
\]

We use in our work the following inequalities for $x\in \left( 0,\pi
/2\right) $ and non-negative integers $m$ and $n,$%
\begin{equation}
u\left( x,2n+1\right) <\sin x<u\left( x,2n\right)   \label{sin}
\end{equation}%
\begin{equation}
v\left( x,2m+1\right) <\cos x<v\left( x,2m\right) ,  \label{cos}
\end{equation}
where%
\[
u\left( x,p\right) =\sum_{k=0}^{p}\frac{\left( -1\right) ^{k}}{\left(
2k+1\right) !}x^{2k+1}\ ,\ \ \ v\left( x,q\right) =\sum_{k=0}^{q}\frac{%
\left( -1\right) ^{k}}{\left( 2k\right) !}x^{2k}
\]%
are truncations of%
\[
\sin x=\sum_{k=0}^{\infty }\frac{\left( -1\right) ^{k}}{\left( 2k+1\right) !}%
x^{2k+1}\ ,\ \ \ \cos x=\sum_{k=0}^{\infty }\frac{\left( -1\right) ^{k}}{%
\left( 2k\right) !}x^{2k}.
\]%
We call a \emph{coefficient-positive polynomial }a\emph{\ }polynomial with
all coefficients positive. For sake of continuity of the presentation, we
provide the involved polynomials in a separate section entitled \emph{%
Appendix }in the final part of this paper.\bigskip 

\emph{Proof of Theorem 1. }We have%
\[
2+\left( \frac{8}{45}-a\left( x\right) \right) x^{3}\tan x-\left( \left( 
\frac{\sin x}{x}\right) ^{2}+\frac{\tan x}{x}\right) <m\left( x\right) ,
\]%
where%
\[
m\left( x\right) =2+\left( \frac{8}{45}-a\left( x\right) \right) x^{3}\frac{%
u\left( x,4\right) }{v\left( x,3\right) }-\left( \left( \frac{u\left(
x,3\right) }{x}\right) ^{2}+\frac{1}{x}\cdot \frac{u\left( x,3\right) }{%
v\left( x,4\right) }\right) 
\]%
Thus%
\[
m\left( \frac{1}{x}\right) =-\frac{A\left( x-1\right) }{76204800x^{12}B%
\left( x-1\right) C\left( x-1\right) },
\]%
where $A,$ $B,$ $C$ are \emph{coefficient-positive. }For every $x\geq 1,$ $%
m\left( 1/x\right) <0,$ that is $m\left( x\right) <0,$ for every $x\in
\left( 0,1\right) .$ The left-hand side inequality (\ref{we}) is proved.

For the other side, we have%
\[
2+\left( \frac{8}{45}-b\left( x\right) \right) x^{3}\tan x-\left( \left( 
\frac{\sin x}{x}\right) ^{2}+\frac{\tan x}{x}\right) >n\left( x\right) ,
\]%
where%
\[
n\left( x\right) =2+\left( \frac{8}{45}-b\left( x\right) \right) x^{3}\frac{%
u\left( x,5\right) }{v\left( x,4\right) }-\left( \left( \frac{u\left(
x,6\right) }{x}\right) ^{2}+\frac{1}{x}\cdot \frac{u\left( x,6\right) }{%
v\left( x,3\right) }\right) 
\]%
Thus%
\[
n\left( \frac{1}{x}\right) =-\frac{D\left( x-1\right) }{38\,\allowbreak
775\,788\,043\,\allowbreak 632\,640\,000x^{24}B\left( x-1\right) C\left(
x-1\right) },
\]%
where $D,$ $B,$ $C$ are \emph{coefficient-positive. }For every $x\geq 1,$ $%
n\left( 1/x\right) <0,$ that is $n\left( x\right) <0,$ for every $x\in
\left( 0,1\right) .$ The right-hand side inequality (\ref{we}) is proved and
the conclusion follows.$\square $\bigskip 

\emph{Proof of Theorem 2. }Since we are interested now for $x$ near $\pi /2$
and estimates (\ref{sin})-(\ref{cos}) give good results near zero, we
replace $x$ by $\pi /2-x$ in the requested inequality to get%
\begin{eqnarray}
&&2+\left( \frac{16}{\pi ^{4}}+c\left( \frac{\pi }{2}-x\right) \right)
\left( \frac{\pi }{2}-x\right) ^{3}\cot x  \label{ss} \\
&<&\left( \frac{\cos x}{\frac{\pi }{2}-x}\right) ^{2}+\frac{\cot x}{\frac{%
\pi }{2}-x}  \nonumber \\
&<&2+\left( \frac{16}{\pi ^{4}}+d\left( \frac{\pi }{2}-x\right) \right)
\left( \frac{\pi }{2}-x\right) ^{3}\cot x  \nonumber
\end{eqnarray}%
For the left-hand side inequality (\ref{ss}), we have%
\[
2+\left( \frac{16}{\pi ^{4}}+c\left( \frac{\pi }{2}-x\right) \right) \left( 
\frac{\pi }{2}-x\right) ^{3}\cot x-\left( \left( \frac{\cos x}{\frac{\pi }{2}%
-x}\right) ^{2}+\frac{\cot x}{\frac{\pi }{2}-x}\right) <p\left( x\right) ,
\]%
where%
\[
p\left( x\right) =2+\left( \frac{16}{\pi ^{4}}+c\left( \frac{\pi }{2}%
-x\right) \right) \left( \frac{\pi }{2}-x\right) ^{3}\frac{v\left(
x,2\right) }{u\left( x,1\right) }-\left( \left( \frac{v\left( x,1\right) }{%
\frac{\pi }{2}-x}\right) ^{2}+\frac{1}{\frac{\pi }{2}-x}\cdot \frac{v\left(
x,1\right) }{u\left( x,2\right) }\right) .
\]%
As%
\[
p\left( \frac{1}{x}\right) =-\frac{E\left( x-3\right) }{2\pi ^{5}x^{5}\left(
\pi x-2\right) ^{2}\left( 6x^{2}-1\right) \left( 120x^{4}-20x^{2}+1\right) },
\]%
with $E$ a \emph{positive-coefficient }polynomial, it results $p\left(
1/x\right) <0,$ for every $x\geq 3.$ Hence $p\left( x\right) <0,$ for every $%
x\in \left( 0,\frac{1}{3}\right) .$ By replacing back $x$ by $\pi /2-x,$ the
left-hand side inequality (\ref{bb}) is proved.

For the right-hand side inequality (\ref{bb}), we have%
\[
2+\left( \frac{16}{\pi ^{4}}+d\left( \frac{\pi }{2}-x\right) \right) \left( 
\frac{\pi }{2}-x\right) ^{3}\cot x-\left( \left( \frac{\cos x}{\frac{\pi }{2}%
-x}\right) ^{2}+\frac{\cot x}{\frac{\pi }{2}-x}\right) >q\left( x\right) ,
\]%
where%
\[
q\left( x\right) =2+\left( \frac{16}{\pi ^{4}}+d\left( \frac{\pi }{2}%
-x\right) \right) \left( \frac{\pi }{2}-x\right) ^{3}\frac{v\left(
x,1\right) }{u\left( x,2\right) }-\left( \left( \frac{v\left( x,2\right) }{%
\frac{\pi }{2}-x}\right) ^{2}+\frac{1}{\frac{\pi }{2}-x}\cdot \frac{v\left(
x,2\right) }{u\left( x,1\right) }\right) .
\]%
As%
\[
q\left( \frac{1}{x}\right) =\frac{F\left( x-2\right) }{144\pi
^{6}x^{6}\left( \pi x-2\right) ^{2}\left( 6x^{2}-1\right) \left(
120x^{4}-20x^{2}+1\right) },
\]%
with $F$ a \emph{positive-coefficient }polynomial, it results $q\left(
1/x\right) >0,$ for every $x\geq 2.$ Hence $q\left( x\right) >0,$ for every $%
x\in \left( 0,\frac{1}{2}\right) .$ By replacing back $x$ by $\pi /2-x,$ the
right-hand side inequality (\ref{bb}) is proved.$\square $

\section{Appendix}

$A\left( x\right) =\allowbreak 22\,755\,\allowbreak
703\,637\,276x+211\,707\,\allowbreak 277\,583\,166\allowbreak
x^{2}+1224\,891\,\allowbreak 661\,206\,384x^{3}+4943\,363\,\allowbreak
812\,472\,396\allowbreak x^{4}+14\,788\,229\,\allowbreak
101\,047\,312x^{5}+34\,004\,439\,\allowbreak 843\,092\,312\allowbreak
x^{6}+61\,480\,986\,\allowbreak 223\,883\,520x^{7}+88\,639\,341\,\allowbreak
361\,781\,280\allowbreak x^{8}+102\,698\,848\,\allowbreak
813\,276\,800x^{9}+95\,874\,412\,\allowbreak 155\,439\,680\allowbreak
x^{10}+71\,962\,038\,\allowbreak
122\,849\,280x^{11}+43\,113\,449\,\allowbreak 986\,871\,040\allowbreak
x^{12}+20\,340\,874\,\allowbreak 557\,158\,400x^{13}+7392\,192\,\allowbreak
479\,884\,800\allowbreak x^{14}+1996\,931\,\allowbreak
543\,040\,000x^{15}+377\,644\,\allowbreak 774\,348\,800\allowbreak
x^{16}+44\,618\,\allowbreak 215\,219\,200x^{17}+2478\,\allowbreak
789\,734\,400\allowbreak x^{18}+1143\,\allowbreak 460\,110\,295$

$B\left( x\right) =\allowbreak
208\,208x+836\,584x^{2}+1861\,440x^{3}+2521\,680\allowbreak
x^{4}+2136\,960x^{5}+1108\,800x^{6}+322\,560x^{7}+\allowbreak
40\,320x^{8}+21\,785$

$C\left( x\right) =2940x+8670x^{2}+12\,960x^{3}+10\,440x^{4}+4320\allowbreak
x^{5}+720x^{6}+389\allowbreak $

$D\left( x\right) =1046\,241\,\allowbreak 057\,892\,766\,\allowbreak
995\,181\,020x+15\,160\,196\,\allowbreak 458\,527\,770\,\allowbreak
766\,797\,390\allowbreak x^{2}+141\,133\,599\,\allowbreak
892\,553\,591\,\allowbreak 885\,361\,888x^{3}+\allowbreak
948\,620\,938\,\allowbreak 108\,112\,884\,\allowbreak
482\,060\,312x^{4}+4904\,616\,140\,\allowbreak 661\,691\,784\,\allowbreak
916\,252\,576\allowbreak x^{5}+20\,\allowbreak 292\,226\,949\,\allowbreak
171\,860\,085\,\allowbreak 229\,532\,656\allowbreak x^{6}+69\,\allowbreak
004\,718\,328\,\allowbreak 773\,496\,670\,\allowbreak
828\,392\,960\allowbreak x^{7}+196\,\allowbreak 572\,398\,646\,\allowbreak
157\,608\,543\,\allowbreak 611\,123\,520\allowbreak x^{8}+475\,\allowbreak
714\,377\,285\,\allowbreak 153\,830\,551\,\allowbreak
505\,007\,360\allowbreak x^{9}+988\,\allowbreak 353\,654\,075\,\allowbreak
203\,455\,545\,\allowbreak 964\,663\,680\allowbreak x^{10}+1776\,\allowbreak
901\,943\,662\,\allowbreak 260\,506\,348\,\allowbreak
259\,502\,080\allowbreak x^{11}+2780\,\allowbreak 857\,159\,731\,\allowbreak
825\,294\,311\,\allowbreak 578\,805\,760\allowbreak x^{12}+3804\,\allowbreak
874\,547\,192\,\allowbreak 905\,665\,416\,\allowbreak
956\,139\,520\allowbreak x^{13}+4565\,\allowbreak 017\,129\,747\,\allowbreak
537\,959\,707\,\allowbreak 485\,025\,280\allowbreak x^{14}+4811\,\allowbreak
255\,903\,914\,\allowbreak 230\,000\,709\,\allowbreak
907\,251\,200\allowbreak x^{15}+4457\,\allowbreak 372\,961\,132\,\allowbreak
476\,471\,420\,\allowbreak 437\,708\,800\allowbreak x^{16}+3628\,\allowbreak
467\,027\,264\,\allowbreak 351\,719\,971\,\allowbreak
654\,041\,600\allowbreak x^{17}+2591\,\allowbreak 354\,540\,244\,\allowbreak
462\,964\,847\,\allowbreak 504\,179\,200\allowbreak x^{18}+1619\,\allowbreak
215\,048\,600\,\allowbreak 627\,699\,570\,\allowbreak
212\,864\,000\allowbreak x^{19}+881\,\allowbreak 619\,928\,396\,\allowbreak
589\,360\,808\,\allowbreak 448\,819\,200\allowbreak x^{20}+415\,\allowbreak
895\,584\,701\,\allowbreak 273\,283\,689\,\allowbreak
578\,496\,000\allowbreak x^{21}+168\,\allowbreak 684\,403\,283\,\allowbreak
843\,009\,965\,\allowbreak 129\,728\,000\allowbreak x^{22}+58\,\allowbreak
222\,827\,644\,\allowbreak 262\,986\,541\,\allowbreak
957\,120\,000\allowbreak x^{23}+16\,\allowbreak 867\,645\,767\,\allowbreak
995\,215\,775\,\allowbreak 662\,080\,000\allowbreak
x^{24}+4025\,191\,509\,\allowbreak 472\,612\,407\,\allowbreak
312\,384\,000\allowbreak x^{25}+770\,502\,378\,\allowbreak
722\,331\,298\,\allowbreak 627\,584\,000\allowbreak
x^{26}+113\,749\,455\,\allowbreak 739\,846\,813\,\allowbreak
286\,400\,000\allowbreak x^{27}+12\,158\,865\,\allowbreak
105\,647\,881\,\allowbreak 420\,800\,000x^{28}+\allowbreak
837\,557\,\allowbreak 021\,742\,465\,\allowbreak
024\,000\,000x^{29}+27\,918\,\allowbreak 567\,391\,415\,\allowbreak
500\,800\,000\allowbreak x^{30}+34\,825\,\allowbreak
460\,317\,038\,\allowbreak 428\,686\,405$

$E\left( x\right) =x^{12}\left( 172\,800\pi ^{4}-17\,280\pi ^{6}\right) $

$+\allowbreak x^{11}\left( 6220\,800\pi ^{4}-921\,600\pi ^{3}+99\,840\pi
^{5}-622\,080\pi ^{6}-960\pi ^{7}\right) +\cdots $

$F\left( x\right) =\left( 5253\,120\pi ^{6}-58\,060\,800\pi ^{4}+69\,120\pi
^{8}\right) x^{12}$

$+\left( 348\,364\,800\pi ^{3}-1393\,459\,200\pi ^{4}-33\,177\,600\pi
^{5}\right. $

$+\left. 126\,074\,880\pi ^{6}-286\,848\pi ^{7}+1658\,880\pi ^{8}\right)
x^{11}+\cdots $

\begin{acknowledgement}
This work was supported by a grant of the Romanian National Authority for
Scientific Research, CNCS-UEFISCDI project number PN-II-ID-PCE-2011-3-0087.
\end{acknowledgement}

\begin{acknowledgement}
Computations in this paper were made using Maple software.
\end{acknowledgement}

\end{document}